\documentclass[12pt,leqno]{article}
\usepackage{amsmath,amssymb,latexsym}
\usepackage{amscd}
\usepackage{mathrsfs}
\usepackage{srcltx}
\usepackage{theorem}
\usepackage{amsfonts}

\newcommand{\N}{\mathbb{N}}        

\newcommand{\To}{\longrightarrow}

\newcommand{\om}{\omega}

\newcommand{\vp}{\varphi}

\newcommand{\lA}{\lambda(A)}

\title{\sc Expos\'e\\
 \large on a paper of Dronov and Kaplitskii
}
\author{Dietmar Vogt}
\date{}

\begin{document}

\maketitle

In \cite{DK} Dronov and Kaplitzki showed that every complemented subspace of a nuclear K\"othe space $E$ with a regular basis of type ($d_1$) has a basis so, in particular, solving the long standing problem whether any complemented subspace of the space (s) of rapidly decreasing sequences has a basis. We
present a slightly modified version of their proof which shows that the range of every closed-range operator in $E$ has a basis.


Let $\lA$ be a nuclear K\"othe space a regular basis of type ($d_1$). The latter means that $E$ has property (DN). Without restriction of generality we may assume:

\begin{enumerate}
\item $a_{1,n}= 1$ for all $n$.
\item $a^2_{k,n}\le a_{k+1,n}$ for all $k,n$.
\item $\sum_{n=1}^\infty \frac{a_{k,n}}{a_{k+1,n}}\le 1$ for all $k$.
\item $\frac{a_{k,n+1}}{a_{k+1,n+1}}\le \frac{a_{k,n}}{a_{k+1,n}}$ for all $k,n$.
\end{enumerate}

Due to nuclearity we may use the following two equivalent norm systems

\begin{itemize}
\item $|x|_k= \sup_n |x_n|a_{k,n}$ , $k\in\N$.
\item $\|x\|_k = \left(\sum_{n=1}^\infty |x_n|^2 a^2_{k,n}\right)^{1/2}$, $k\in\N$.
\end{itemize}

Because of 3. above we obtain:
\begin{itemize}
\item $|x|_k\le \|x\|_k\le |x|_{k+1}$, $k\in\N$.
\end{itemize}

The respective local Banach spaces are
$$G_k=c_0(a_k)\text{ and }H_k=\ell_2(a_k),$$

We consider an operator $T\in L(E)$ and we set $F=T(E)\subset E$, We want to study properties of $F$. $T$ is given in the form
$$Tx=\left( \sum_{j=1}^\infty t_{i,j} x_j\right)_{i\in\N}.$$
We define an operator $|T|$ by
$$|T|x=\left( \sum_{j=1}^\infty |t_{i,j}| x_j\right)_{i\in\N}.$$
To see that this defines an operator $|T|\in L(E)$ we recall the explicit description of the matrices of operators in $c_0$. $(t_{i,j})_{i,j\in\N}$ defines an operator in $c_0$ iff 1. $\sup_i \sum_{j=1}^\infty |t_{i,j}|<\infty$ and 2. $\lim_i t_{i,j}=0$ for all $j$. This description depends on $|t_{i,j}|$ only.

Without restriction of generality we may assume
\begin{enumerate}
\item[5.\,] $\| |T|x\||_k\le \frac{1}{2} |x|_{k+1}$, $k\in\N$.
\end{enumerate}

Next we define three ``dead-end spaces'' that is continuously imbedded Banach spaces, given by two weights 
$$a_{\infty,n}^2:= \sum_{k=1}^\infty \delta_k^2 a_{k,n}^2,$$
where the $\delta_k$ will be determined later, and
$$b_{\infty,n}:=a_{n,n}.$$
We set $H_\infty=\ell_2(a_\infty)$ with the norm
$$\|x\|_\infty^2:=\|x\|_{H_\infty}^2= \sum_{k=1}^\infty \delta_k^2 \|x\|_k^2.$$
We set $G_\infty=c_0(a_\infty)$ with the norm
$$|x|_\infty:=|x|_{G_\infty}= \sup_n |x_n|a_{\infty,n}.$$
Moreover we set $G_{\infty,0}=c_0(b_\infty)$ with the norm
$$|x|_{\infty,0}:=|x|_{G_{\infty,0}}=\sup_n |x_n| b_{\infty,n}.$$
We obtain
$$G_{\infty,0}\subset H_k \subset G_k, \qquad k\in\N.$$
The second inclusion is obvious, the first one we get from
$$\sum_{n=k+1}^\infty a_{k,n}^2 |x_n|^2 = \sum_{n=k+1}^\infty\frac{a_{k,n}^2}{a_{n,n}^2}\, a_{n,n}^2|x_n|^2  \le \left(\sum_{n=k+1}^\infty\frac{a_{k,n}^2}{a_{n,n}^2}\right) |x|^2_{\infty,0}\le |x|^2_{\infty,0}.$$

With some constant $D_k$ we have
$$\|x\|_k\le D_k |x|_{\infty,0}.$$
We may assume $D_k\le D_{k+1}$ for all $k$. We obtain
\begin{eqnarray*}
\|Tx\|_{\infty}^2 &=& \sum_{k=1}^\infty \delta_k^2 \|Tx\|_k^2
\le \sum_{k=1}^\infty \delta_k^2  \|x\|_{k+1}^2
\le \left(\sum_{k=1}^\infty \delta_k^2 D_{k+1}^2\right) |x|_{\infty,0}^2
\le |x|_{\infty,0}^2,
\end{eqnarray*}
where we have chosen, also under consideration of later application, $\delta_k\le 1/(2^k D_{k+2})$. So we have shown
\begin{equation}\label{e1}
\|Tx\|_{\infty}\le |x|_{\infty,0}
\end{equation}

This shows that $L:=\{Tx\,:\,x\in G_{\infty,0}\}\subset H_\infty\}$, We set $F_k$ the completion of $L$ with respect to $\|\cdot\|_k$ and $F_\infty$ the completion with respect to $\|\cdot\|_\infty$. The embedding $J:F_\infty\hookrightarrow F_1$ is clearly nuclear, hence we can expand it as
$$J(x) = \sum_{j=1}^\infty \langle x,f_j\rangle_{H_1} f_j$$
where $(f_j)_{j\in\N}$ is orthogonal in $H_\infty$, orthonormal in $H_1$. We set
$$T_n(x) = \sum_{j=1}^n \langle Tx,f_j\rangle_{H_1} f_j.$$
For every $x\in G_{\infty,0}$ we get $T_n(x)\To T(x)$ in $H_\infty$. We want to show that the family of maps $\{T_n\,:\,n\in\N\}$ is equicontinuos in $L(E)$. This will imply
\begin{equation*}
Tx=\sum_{j=1}^\infty \langle Tx,f_j\rangle_{H_1} f_j
\end{equation*}
for all $x\in E$, the series converging in $E$.

\bf Theorem: \it If $T(E)$ is closed, then $F_\infty\subset T(E)$, hence all $f_j\in T(E)$. So the $f_j$ are a basis in $T(E)$. \rm

Because of orthogonality we have
\begin{itemize}
\item $\|T_n x\|_1\le \|Tx\|_1$.
\item $\|T_n x\|_\infty \le \|Tx\|_\infty$.
\end{itemize}
For $x\in\vp^+$ we have $\|Tx\|_1\le \| |T|x\|_1$ and therefore we have
\begin{itemize}
\item $\|T_nx\|_1\le \| |T|x\|_1$.
\end{itemize}

To get an estimate between the $\|\cdot\|_\infty$ and the $|\cdot|_\infty$ norm, we fix some $r\in\N$ and obtain
\begin{eqnarray}\label{e2}
\|x\|_\infty&=&\|x\|_{\ell_2(a_\infty)}\le \|x\|_{\ell_1(a_\infty)}=\sum_{n=1}^\infty |x_n| a_{\infty,n}
=  \sum_{n=1}^\infty \frac{a_{r,n}}{a_{r+1,n}}\left( \frac{a_{r+1,n}}{a_{r,n}} |x_n| a_{\infty,n}\right)\\ \nonumber
&\le& \sum_{n=1}^\infty \frac{a_{r,n}}{a_{r+1,n}}\sup_{n\in\N}\left\{ \frac{a_{r+1,n}}{a_{r,n}} |x_n| a_{\infty,n}\right\}\\ \nonumber
&\le& \left|\frac{a_{r+1}}{a_r}x\right|_\infty \nonumber
\end{eqnarray}
for all $x\in G_{\infty,0}$.

To see that $\frac{a_{r+1}}{a_r}x\in G_\infty$ for all $x\in G_{\infty,0}$ we use the estimate:
$$\|\frac{a_{r+1}}{a_r}x\|_k^2\le \|a_k x\|_k^2=\sum_{n=1}^\infty a_{k,n}^4 |x_n|^2 \le \sum_{n=1}^\infty a_{k+1,n}^2 |x_n|^2 =\|x\|_{k+1}^2$$
for all $k>r$. From that we obtain
\begin{eqnarray*}
\|\frac{a_{r+1}}{a_r}x\|_\infty^2 &=& \sum_{k=1}^\infty \delta_k^2 \|\frac{a_{r+1}}{a_r}x\|_k^2\\
&\le& \sum_{k=1}^r \delta_k^2 \|\frac{a_{r+1}}{a_r}x\|_k^2+\sum_{k=r+1}^\infty \delta_k^2 \|\frac{a_{r+1}}{a_r}x\|_k^2\\
&\le& \|\frac{a_{r+1}}{a_r}x\|_r^2 \sum_{k=1}^r \delta_k^2 + \sum_{k=r+1}^\infty\delta_k^2 \|x\|_{k+1}^2\\
&=&  \|x\|_{r+1}^2 \sum_{k=1}^r \delta_k^2 + \sum_{k=r+1}^\infty\delta_k^2 \|x\|_{k+1}^2\\
&\le& (D_{r+1}^2+\sum_{k=r+1}^\infty\delta_k^2 D_{k+1}^2) |x|_{\infty,0}^2\\
&\le& (D_{r+1}^2+1) |x|_{\infty,0}^2
\end{eqnarray*}
Therefore $\frac{a_{r+1}}{a_r}x\in H_\infty\subset G_\infty$ for all $x\in G_{\infty,0}$.

We apply the previous to $|T|x$ and obtain:
\begin{eqnarray*}
\|\frac{a_{r+1}}{a_r}|T|x\|_\infty^2 &\le&\||T|x\|_{r+1}^2 \sum_{k=1}^r \delta_k^2 + \sum_{k=r+1}^\infty\delta_k^2 \||T|x\|_{k+1}^2\\
&\le&\|x\|_{r+2}^2 \sum_{k=1}^r \delta_k^2 + \sum_{k=r+1}^\infty\delta_k^2 \|x\|_{k+2}^2\\
&\le& (D_{r+2}^2+\sum_{k=r+1}^\infty\delta_k^2 D_{k+2}^2) |x|_{\infty,0}^2\\
&\le& (D_{r+2}^2+1) |x|_{\infty,0}^2.
\end{eqnarray*}

We define $J_r x:= \frac{a_{r+1}}{a_r}x$ and we have shown
\begin{equation}\label{e4}
\|J_r |T|x\|_\infty\le M_r |x|_{\infty,0}
\end{equation}
with $M_r^2 = D_{r+2}^2+1$.

We will use the canonical projection in the sequence space $E$:
$$\left(Q^{(N)}x\right)_n =\begin{cases} x_n, & n=1,2,\dots, N \\
0, & n>N.
\end{cases}$$

We obtain for all $x\in G_1$
\begin{eqnarray}\label{e3}
\left|\frac{a_{r+1}}{a_r}\,|T|\, \frac{a_1}{a_{r+2}}\, Q^{(N)} x\right|_r &=&
\left|\,|T|\, \frac{a_1}{a_{r+2}}\, Q^{(N)} x\right|_{r+1}\\ \nonumber
&\le& \frac{1}{2} \left|\frac{a_1}{a_{r+2}}\, Q^{(N)} x\right|_{r+2} = \frac{1}{2}|Q^{(N)} x|_1\le\frac{1}{2}|x|_1.
\end{eqnarray}

We set $J_r'= a_1/a_{r+2}$ and define:
$$A_r=J_r |T| J_r' \quad\text{ and }\quad A_r^{(N)} = A_r\circ Q^{(N)}.$$
The operator acts $G_1\to G_r$, hence also $G_r\to G_r$ and $\|A_r^{(N)}\|_{G_r\to G_r}\le \frac{1}{2}$.

For all $x\in G_{\infty,0}^+ = \om^+\cap G_{\infty,0}$  we obtain
$$|T_nx|_\infty\le \|T_n x\|_\infty\le \|Tx\|_\infty \le \|\,|T|x\|_\infty \le \left|\frac{a_{r+1}}{a_r}\,|T|\,x\right|_\infty= |J_r |T|x|_\infty.$$
The last inequality comes from (\ref{e2}).

For $x\in\vp^+$ we get the estimates:
$$|T_nJ_r'x|_\infty \le |J_r |T| J_r' x|_\infty =|A_rx|_\infty$$
$$|T_n J_r'x|_1 \le |J_r |T| J_r' x|_1\le |x|_1 $$

In the second line the first estimate comes from
$$\|x\|_1^2=\sum_n|x_n|^2 = \sum_n \frac{a^2_r}{a^2_{r+1}} \frac{a^2_{r+1}}{a^2_{r}} |x_n|^2\le |J_r x|_1^2$$
which implies
$$|T_n x|_1\le \|T_n x\|_1\le \|Tx\|_1\le \|\,|T|x\|_1\le |J_r |T| x|_1.$$
The second estimate is (\ref{e3}).

We define
$$S_{n,r}=T_n J_r' \quad \text{and }\quad S_{n,r}^{(N)}= S_{n,r}Q^{(N)}.$$

The first of the inequalities above is not applicable for interpolation, but it becomes applicable when we restrict it to a suitable cone.

We set
$$Q_{r,N}=\{x\in \om^+\,:\, x\ge A_r^{(N)}x\}.$$
Then $Q_{r,N}$ is a cone and we have:
$$|S_{n,r}^{(N)}x|_1 \le |x|_1 \text{ for }x\in G_1^+,$$
$$|S_{n,r}^{(N)}x|_\infty\le |x|_\infty \text{ for } x\in Q_{r,N}\cap G_\infty^+.$$

By use of the interpolation theorem for cones \cite[Theorem 1]{DK} we obtain that for every $r$ there is a constant $C(r)$ such that
$$|S_{n,r}^{(N)}x|_r\le C(r)\, |x|_r \text{ for } x\in Q_{r,N}\cap G_r^+.$$

Since $A_r^{(N)}\in L(G_r)$ with $\|A_r^{(N)}\|\le 1/2$ the operator $B:=I-A_r^{(N)}$ is invertible in $L(G_r)$ and $\|B^{-1}\|\le 2$. We can write for $x\in G_r$
$$x= B^{-1}Bx=B^{-1}(Bx)_+ - B^{-1}(Bx)_-.$$
Clearly $B^{-1}(Bx)_+,B^{-1}(Bx)_-\in Q_{r,N}\cap G_r^+$. Positivity is seen by the Neumann series.
Since $\|B^{-1}B\|\le 4$, we have shown that every $x\in G_r$ can be written as $x=y-z$ here $y,z\in Q_{r,N}\cap G_r^+$ and $\|y\|_r\le 4\|x\|_r$ and $\|z\|_r\le 4\|x\|_r$.

It follows that
$$|S_{n,r}^{(N)}x|_r\le 8C(r)\, |x|_r \text{ for } x\in G_r.$$
With $N\to \infty$ we conclude that
$$|T_nJ_r'x|_r \le 8C(r)\, |x|_r \text{ for } x\in G_r.$$
Finally we have
$$|T_n x|_r \le 8C(r)\,|a_{r+2} x|_r \le 8C(r)|x|_{r+3}.$$
Since that can be done for all $r$ the family $(T_n)_{n\in\N}$ is equicontinuous in $E$ which had to be shown.

We have to verify the conditions of the cone interpolation theorem \cite[Theorem 1]{DK}.

1. $Q_{r,N}$ is a lower semi-lattice. This follows immediately from the fact that $A_r^{(N)}$ is monotone on $\om^+$.

2. $Q_{r,N}\cap G_\infty^+$ is total in $G_\infty$. We have seen that $I-A_r^{(N)}$ is invertible in $G_r$, hence injective on $G_\infty$. We show that $A_r^{(N)}\in L(G_\infty)$. Then, as $A_r^{(N)}$ is finite dimensional, $I-A_r^{(N)}$ is bijective in $G_\infty$. Arguing as above we obtain that $Q_{r,N}\cap G_\infty^+-Q_{r,N}\cap G_\infty^+=G_\infty$.

We have the following chain of inequalities, where the second one is (\ref{e4}) the last one finite dimensionality of $Q^{(N)}$.
$$|A_r^{(N)}x|_\infty=|J_r |P| J_r' Q^{(N)}x|_\infty \le \|J_r |P| J_r' Q^{(N)}x\|_\infty \le
|J_r' Q^{(N)}x|_{\infty,0} \le | Q^{(N)}x|_{\infty,0}\le C(N) |x|_\infty.$$
This completes the argument for 2.

3. Finally we have to show that $Q_{r,N}\cap G_\infty^+$ contains a strictly positive element. To show that we choose a strictly positive element $x_0\in G_\infty$ and put $x=(I-A_r^{(N)})^{-1} x_0$. $x$ can be calculated by means of the Neumann series in $G_r$. Since $A_r^{(N)}$ is positive this shows that $x\ge x_0$ which shows the claim.


\vspace{.5cm}

\noindent Bergische Universit\"{a}t Wuppertal,
\newline FB Math.-Nat., Gau\ss -Str. 20,
\newline D-42119 Wuppertal, Germany
\newline e-mail: dvogt@math.uni-wuppertal.de

\end{document}